\documentclass[12pt]{iopart}
\usepackage{verbatim,color,amssymb,epsfig, multirow,iopams,amsfonts,cite,array,graphicx}
\usepackage{subfigure,float}
\graphicspath{{./Figs/}}

\def\S{{\cal D}}

\def\boxit#1{\vbox{\hrule\hbox{\vrule\kern6pt
          \vbox{\kern6pt#1\kern6pt}\kern6pt\vrule}\hrule}}

\def\log{\hbox{log}}

\def\JRSSB{{\it Journal of the Royal Statistical Society, Series B}}

\def\JASA{{\it Journal of the American Statistical Association}}

\def\JASA{{\it Journal of the American Statistical Association}}

\def\JRSSB{{\it Journal of the Royal Statistical Society, Series B}}

\def\JRSSB{{\it Journal of the Royal Statistical Society, Series B}}

\def\JASA{{\it Journal of the American Statistical Association}}

\def\bse{\begin{eqnarray*}}
\def\ese{\end{eqnarray*}}
\def\be{\begin{eqnarray}}
\def\ee{\end{eqnarray}}
\def\bq{\begin{equation}}
\def\eq{\end{equation}}
\def\bse{\begin{eqnarray*}}
\def\ese{\end{eqnarray*}}
\def\pr{\hbox{pr}}

\def\th{^{th}}

\begin{document}
\title[Detecting a small emission source]{A Bayesian approach to detection of small low emission sources}
%

\author{Xiaolei Xun$^1$, Bani Mallick$^1$, Raymond J. Carroll$^1$, Peter Kuchment$^2$}
\address{$^1$ Department of Statistics, Texas A\&M University, 3143 TAMU, College Station, TX 77843-3143}
\address{$^2$ Department of Mathematics, Texas A\&M University, 3368 TAMU, College Station, TX 77843-3368}
\ead{\mailto{xxun@stat.tamu.edu}, \mailto{bmallick@stat.tamu.edu}, \mailto{carroll@stat.tamu.edu}, \mailto{kuchment@math.tamu.edu} }

\begin{abstract}
The article addresses the problem of detecting presence and location of a small low emission source inside of an object, when the background noise dominates. This problem arises, for instance, in some homeland security applications.
The goal is to reach the signal-to-noise ratio (SNR) levels on the order of $10^{-3}$.
A Bayesian approach to this problem is implemented in 2D.
The method allows inference not only about the existence of the source,
but also about its location. We derive Bayes factors for model selection and
estimation of location based on Markov Chain Monte Carlo (MCMC) simulation.
A simulation study shows that with
sufficiently high total emission level,
our method can effectively locate the source.
\end{abstract}
\submitto{\IP}

\footnote{ \baselineskip=10pt The correspondence author in the summer will be P.Kuchemnt. If it is accepted, then the correspondence author will be X.Xun.}.

\maketitle


\section{Introduction} \label{sec:intro}

We consider the problem of detecting existence of a low emission radiating source inside a volume, in the presence of a strong random background. One can easily imagine possible applications of such detection, for instance to homeland security. We are interested in the situation when about $99.9\%$ of the total detections come from the background particles\footnote{Most of the background particles are not emitted inside the object, but rather are present in the surrounding environment (e.g., cosmic rays).} and from the particles emitted by the source that have been scattered (we will consider the latter as a part of the background). In other words, only about $1\%$ of detected hits are by the ballistic particles coming from the source. Although medical emission tomographic imaging faces similar problems (e.g. \cite{Herman}), the overwhelming level of noise that has just been mentioned would be considered impossible to handle there. So, how can one attack this problem? Although there is probably no general solution, in the applications we have in mind, the radiating source, if present, would be significantly (on the order of hundred times) geometrically smaller than the whole object. As is explained in \cite{ADHK}, this, and the availability of detectors determining direction of an incoming particle, make detection conceivable under appropriate conditions. In this text, we consider the $2D$ problem. Unlike \cite{ADHK}, where more analytic techniques are considered, we propose a Bayesian method, which allows inference about the existence and location of the possible source.

\begin{figure}[t]
\caption{\label{fig:fig1} 
A direction sensitive (e.g., collimated or Compton type) detector determines the normal parameters $(\theta, S)$ of the trajectory of the incoming particle. The detected particles might be either emitted from the source or coming from random background.}
\begin{indented}
\item[] \scalebox{0.6}{\includegraphics{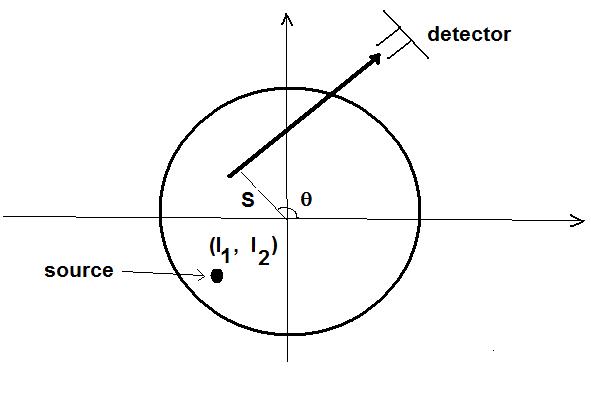}}
\end{indented}

\end{figure}

The problem can be stated as follows. One is interested in certain type of particles (say, $\gamma$-photons or neutrons, although the type of particles is irrelevant for our purpose). Suppose that the observed area belongs to the unit disk D in the $l_1$-$l_2$-plane (see Figure \ref{fig:fig1}). Detectors, placed around the object, are assumed to be able to determine the linear trajectory of each incoming particle. It is assumed that detectors surround the object in such a way that any escaping particle hits a detector (this assumption could in principle be weakened). Most of (or maybe all) detected particles are coming from a random background (and, in particular, are not emitted inside the object). Besides the background emission, a small (in comparison with the total object's size) and low emission source might be present. Many of the particles emitted by the source will be scattered, and only a small number of them will reach the detectors unscattered (ballistic). The goal is to detect the presence of such an object, if the emission is dominated by the background, e.g., such that the ballistic particles coming from a possible source could account to about $0.1\%$ of the total emission.

The set-up is illustrated in Figure \ref{fig:fig1}. The trajectory of a particle that hits the detector can be identified by its normal coordinates $(\theta, S)$, and thus we assume that the detector provides the values $(\theta,S)$ for each hit.
We expect that a radiation source of a small radius $d$ might be present,
in which case we denote its location point as $L=(\ell_1, \ell_2)$.
If a particle is emitted from this source and reaches the detectors ballistically (unscattered), then $\theta$ and $S$ satisfy the inequality
\begin{equation}\label{eq:through}
|\ell_1\cos(\theta)+\ell_2\sin(\theta)-S| \leq d.
\end{equation}
Most particles from the random background normally will not satisfy this condition, but a small portion might. The idea is that if a source is present, ballistic particles coming from it might lead to a statistically significant increase in the number of trajectories satisfying (\ref{eq:through}), and thus to detection. Under appropriate conditions (geometrically sufficiently small source and sufficiently large total count in the sample), this happens to be the case (see the discussion in \cite{ADHK}).

An outline of the paper is as follows.
Section \ref{sec:model} introduces the candidate models
describing the situations without a source and with a source, respectively.
Section \ref{sec:BF} explains the calculation of
Bayes factors for determining the presence of a source, as well as the computational details of our Markov Chain Monte Carlo (MCMC) algorithm. Readers interested mostly in the viability of the approach can skip the subsections, going directly to Section \ref{sec:simu} with the results of a simulation study, where the method is examined for various levels of source emission rate. These results confirm the possibility of detection.
The article ends with conclusions and remarks.

\section{Models} \label{sec:model}

Suppose that the direction sensitive detectors registered hits by $n$ particles and recorded the corresponding normal coordinates
$(\theta_i, S_i)$ for $i=1,\cdots,n$ of their incoming directions.
We denote by $\delta_i$ the (unobserved) indicator that
the $i\th$ particle is coming ballistically from the suspected source.  In other words, $\delta_i=1$ if the $i\th$ particle comes from the source. Otherwise, $\delta_i=0$.

If there is no source present, then $\pr(\delta_i=1)=0$.
If there is a source, we assume that the $\delta_i$'s are independently and
identically distributed according to the $\hbox{Bernoulli}(p)$ law.
This covers also the possible absence of a source, in which case $p=0$.

Our plan is to develop statistical models for each of these cases,
and then decide, based on the collected data and the value of the corresponding
Bayes factor, which model fits better the collected data.

We would like to point out that this is not expected to be a simple problem, since it involves inference as to whether a non-negative parameter $p$ takes its boundary value $p=0$. Even in simple variance components models, frequentist boundary value testing is a difficult matter, see for example \cite{CrRu}.

\subsection{The Model without a Source}\label{subsec:reduced}

When there is no source (we call this model $M_1$), all hits at the detectors come from the random background and thus $\delta_i =0$ for all $i=1,\dots,n$.
We will assume in this text that the random background is isotropic and uniform. In other words, the angle $\theta$ and the distance $S$ from the origin of the trajectory are uniformly distributed:
\begin{equation}\label{eq:null}
\eqalign{
   [ \theta_{i}|\delta_i=0 ] = {\rm Uniform}(0, 2\pi); \cr
   [ S_{i}|\delta_i=0] = {\rm Uniform}(-1, 1).}
   \end{equation}
Notice that particles having trajectories with $|S|>1$ do not get detected and thus do not enter the model.

\subsection{The Model with a Source}\label{subsec:full}

If a source exists (model $M_2$), then $\pr (\delta_i=1) = p > 0 $.
If the particle comes from the background, then $\delta_i=0$ and
relations (\ref{eq:null}) still hold. Assuming that the source in question is isotropic and uniform, when $\delta_i=1$, we have the following distributions:
\begin{equation}
\eqalign{
   [ \theta_{i} | \delta_i=1  ]
         = {\rm Uniform}(0, 2\pi);     \cr
   [ S_{i} | \theta_{i}, (\ell_1, \ell_2), \delta_i=1 ]
         = \ell_1\cos(\theta_{i}) + \ell_2\sin(\theta_{i})
            + {\rm Uniform}(-d,d). }
   \label{eq:full}
\end{equation}
Our goal thus is to choose between the models $M_1$ and $M_2$, based on the measured data.
We show in the following subsection the priors, likelihood and the posterior distribution associated with model $M_2$.

\subsubsection{Likelihood and Posterior of Model with a Source}\label{subsubsec:full}

For the model $M_2$ with a source, see (\ref{eq:full}), parameters of interest are
$\phi = (p, \ell_1, \ell_2)$. Priors are assigned as follows:
\begin{equation*}
 L = (\ell_1, \ell_2) =  \hbox{Uniform} \{\ell_1^2+\ell_2^2 \leq 1\},
\end{equation*}
i.e., \textit{a priori} the source can be located anywhere with equal probability. Also, we allow
$p$ to have the uniform discrete distribution on the set $ \{0<p_1, \cdots, p_h \} $,
which contains $h$ values equally spaced on the interval $[p_1, p_h]$. This set is chosen based on \textit{a priori} estimates of the possible emission strength of the source.

For the implementation of the MCMC algorithm, it is convenient (see Section \ref{subsec:gibbs}) to use the polar coordinates $(r, u)$ of the center $L$:
\begin{equation*}
\ell_1=r\cos u,\quad \ell_2=r\sin u.
\end{equation*}
Denoting the new parameterization by $\psi=(p, r, u)$, the prior $f(\psi|M_2)$ is
\begin{equation*}
f(\psi|M_2)=f(p,r,u|M_2)=h^{-1}\pi^{-1} r\mbox{I}(0\leq r \leq 1) \mbox{I}(0\leq u \leq 2\pi),
\end{equation*}
where $\mbox{I}(\cdot)$ is the indicator function.

Let $Y_i=(\theta_i, S_i)$ denote the $i\th$ observation (i.e., the normal coordinates of the trajectory (at the detector) of the $i\th$ detected particle) and, as before,
let $\delta_i$ be the (unobserved) indicator associated with it.
It will be convenient to introduce vectors $\widetilde{Y}=(Y_1, \dots, Y_n)$ and $\widetilde{\delta} = (\delta_1, \dots, \delta_n)$. Then the likelihood function is
\begin{eqnarray}
 &  f(\widetilde{Y}|p, r, u, M_2 )   \label{eq:large} \\
 = & \hbox{$\prod\limits_{i=1}^n$}
              (4\pi)^{-1} \left[ pd^{-1} \mbox{I} \left\{
                |  r\cos(u)\cos(\theta_i)+r\sin(u)\sin(\theta_i)-S_i  | \leq d \right\} + 1-p  \right]  \nonumber \\
         = & (4\pi)^{-n}   \hbox{$\prod\limits_{i=1}^n$}
                \left[ pd^{-1} \mbox{I} \left\{
                |  r\cos(u-\theta_i)-S_i  | \leq d \right\} + 1-p  \right] \nonumber \\
         = & (4\pi)^{-n}
                (1-p)^{n-J} ( pd^{-1}+1-p ) ^J,\nonumber
\end{eqnarray}
where $J=\hbox{$\sum_{i=1}^n$} \mbox{I} \left\{| r\cos(u-\theta_i)-S_i|\leq d \right\}$ counts the total number of particles whose incoming trajectories at the detectors pass near the location $(r, u)$.

Given the above prior and likelihood, the posterior is
\begin{equation}
 f(\psi|\widetilde{Y}, M_2) = f(p, r, u | \widetilde{Y}, M_2) \propto r(1-p)^{n-J} ( pd^{-1}+1-p ) ^J,
\label{eq:post}
\end{equation}
where $p \in \{p_1, \cdots, p_h \}$, $r\in [0,1]$ and $u \in [0, 2\pi]$.

\section{Model Selection via Bayes Factors} \label{sec:BF}

Let $\pr(M_j)$ be the prior probability of model $M_j$ and $\pr(\widetilde{Y}|M_j)$ be the marginal distribution of the data, given model $M_j$, where $j=1, 2$. We also denote by $\pr(M_j|\widetilde{Y})$ the posterior probability of the model $M_j$.

The parameters of interest are $p$ and $ L=( \ell_1, \ell_2)$. Indeed, if $p=0$, then there is no source, while if $p>0$, the source is present at the location $L$.

In the next section, we describe the Bayesian approach that will be used for model selection. Then the computation and algorithm will be explained.

We use a Bayes factor approach to select between the two models $M_1$ and $M_2$ in question.
The Bayes factor is defined as the ratio of the prior and posterior odds:
\bse
BF & = & \frac{\pr(M_1)/\pr(M_2)}{\pr(M_1|\widetilde{Y})/\pr(M_2|\widetilde{Y})} \\
   & = & \frac{\pr(M_1)\pr(M_2)\pr(\widetilde{Y}|M_2)}{\pr(M_2)\pr(M_1)\pr(\widetilde{Y}|M_1)}
     =   \frac{\pr(\widetilde{Y}|M_2)}{\pr(\widetilde{Y}|M_1)}.   \nonumber
\ese
This number serves as an indicator of which of the models $M_1$ and $M_2$ is more supported by the data.
If $BF>1$, this indicates $M_2$ being more strongly supported by the data. Otherwise, $M_1$ is more strongly supported.
Furthermore, the magnitude of the Bayes factor is a measure of how strong the evidence is for or against $M_1$. According to Kass and Raftery \cite{KaRa}, when the Bayes factor exceeds $3$, $20$ and $150$, one can say that, correspondingly, a positive, strong, and overwhelming evidence exists that a source is present. See \cite[Appendix B]{Je}, \cite{Ev} and \cite{Go} for further interpretation of Bayes factors.

We thus need the marginal distributions $\pr(\widetilde{Y}|M_j)$ to be calculated for each candidate model $M_j, j=1,2$.

Under the null model $M_1$, in which there is no source, one concludes that the corresponding marginal probability density of $\widetilde{Y}$ is:
\begin{equation*}
\pr(\widetilde{Y}|M_1) = (4\pi)^{-n}.
\end{equation*}
When there is a source, we denote by $\Psi$ the {\bf sample space of parameters under} $M_2$. The points $\psi$ from this space is the triples $\psi=(p, r, u)$ described in Section \ref{subsubsec:full}.
Then the marginal probability of $\widetilde{Y}$ under $M_2$ is
\begin{equation*}
\pr(\widetilde{Y}|M_2) = \int\limits_\Psi \pr (\widetilde{Y}, \psi |M_2) f(\psi|M_2) d\psi,
\end{equation*}
a quantity that cannot be computed explicitly.

The usual Monte Carlo method of computation is as follows.
Suppose we have $k=1, \cdots, K$ samples $\psi^{(k)} = (p^{(k)}, r^{(k)}, u^{(k)} )$
from the posterior distribution.
The marginal distribution $\pr(\widetilde{Y}|M_2)$ can then be estimated as
\bse
\widehat{\pr} (\widetilde{Y} | M_2) \approx \{ K^{-1} \sum_{k=1}^K
     \pr (\widetilde{Y}|\psi^{(k)},M_2) ^{-1} \} ^{-1},
\ese
i.e., {\bf the harmonic mean of the likelihoods} $\pr (\widetilde{Y}|\psi, M_2)$ (see, e.g. \cite{KaRa}).
Given this, the Bayes factor is calculated as
\be
    \widehat{\hbox{BF}}
        &=& \widehat{\pr}(\widetilde{Y}|M_2) / \hbox{pr} (\widetilde{Y}|M_1)  \label{eq:BFcalc}  \\
        &=& \left[ K^{-1} \hbox{$\sum_{k=1}^K$}
            \left\{ (1-p^{(k)})^{n-J(k)} (p^{(k)}d^{-1}+1-p^{(k)})^{J(k)}
        \right\}^{-1}    \right]^{-1} ,   \nonumber
\ee
where $J(k)=\hbox{$\sum_{i=1}^n$} \mbox{I}  \left(|r^{(k)} \cos(u^{(k)}-\theta_i)-S_i|\leq d  \right) $. See \cite{KaRa} and \cite{Ra} for more details about calculation of Bayes factors.

Computational details of the MCMC implementations and
calculation of Bayes factors are shown in the following subsections.

\subsection{Computation}\label{subsec:gibbs}

With the model and prior in Section \ref{subsec:full}, the posterior distribution is not straightforward to sample from.
Thus, the Markov Chain Monte Carlo method is used to simulate the parameters from the posterior.

\begin{figure}[t]
\caption{\label{fig:multimode} 
Snapshots of $f(\ell_1, \ell_2 | p, \widetilde{Y})$. The data set contains $200,000$ samples. The true emission rate is p = 0.001, and the source is located at (0.3, 0.6). The left figure is conditioned at $p=0.0002$; the right one assumes $p=0.001$. The multimodality illustrates the difficulty of MCMC sampling in this problem.}
\begin{indented}
\item[] \scalebox{0.3}{\includegraphics{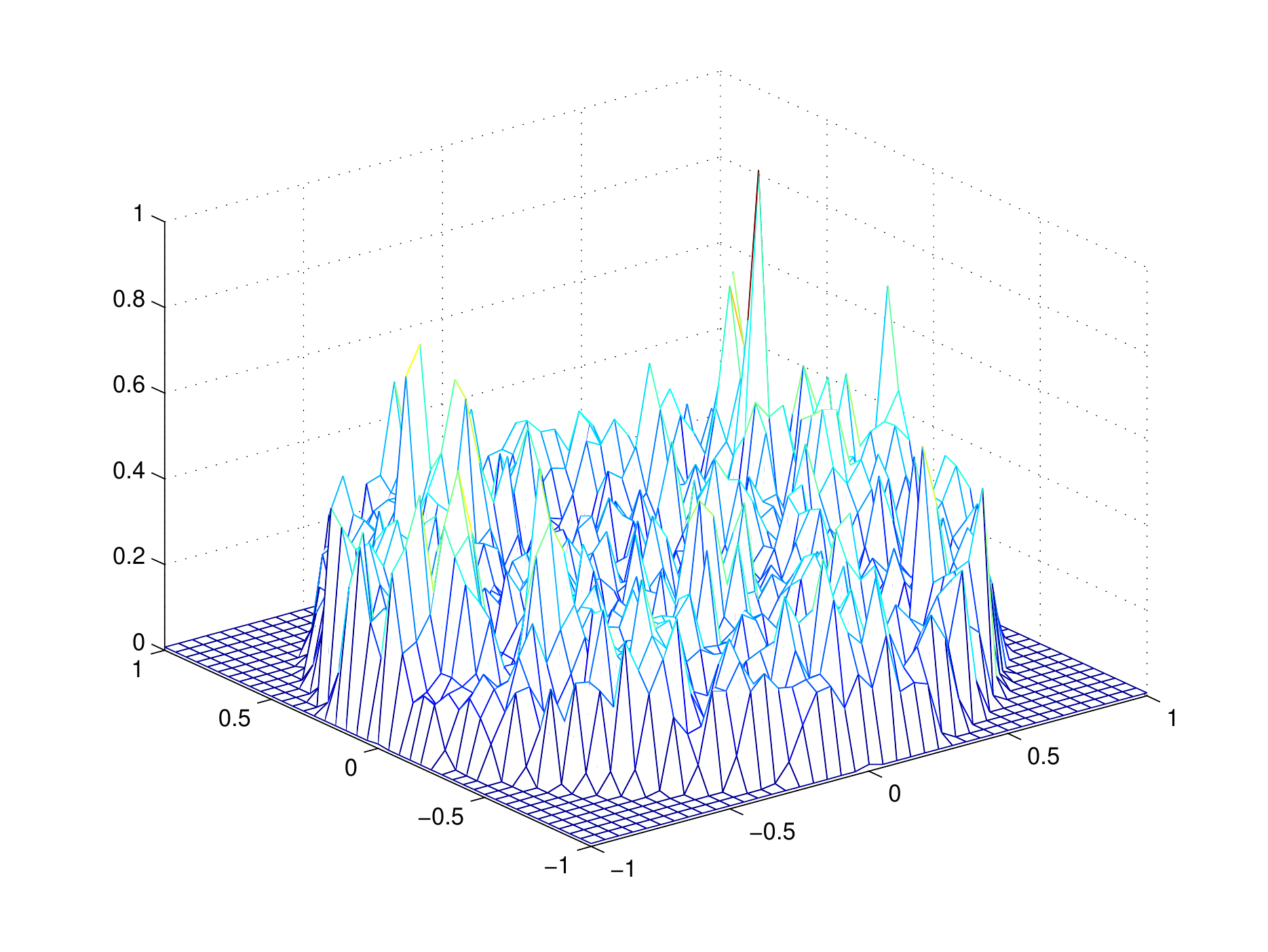}
\includegraphics{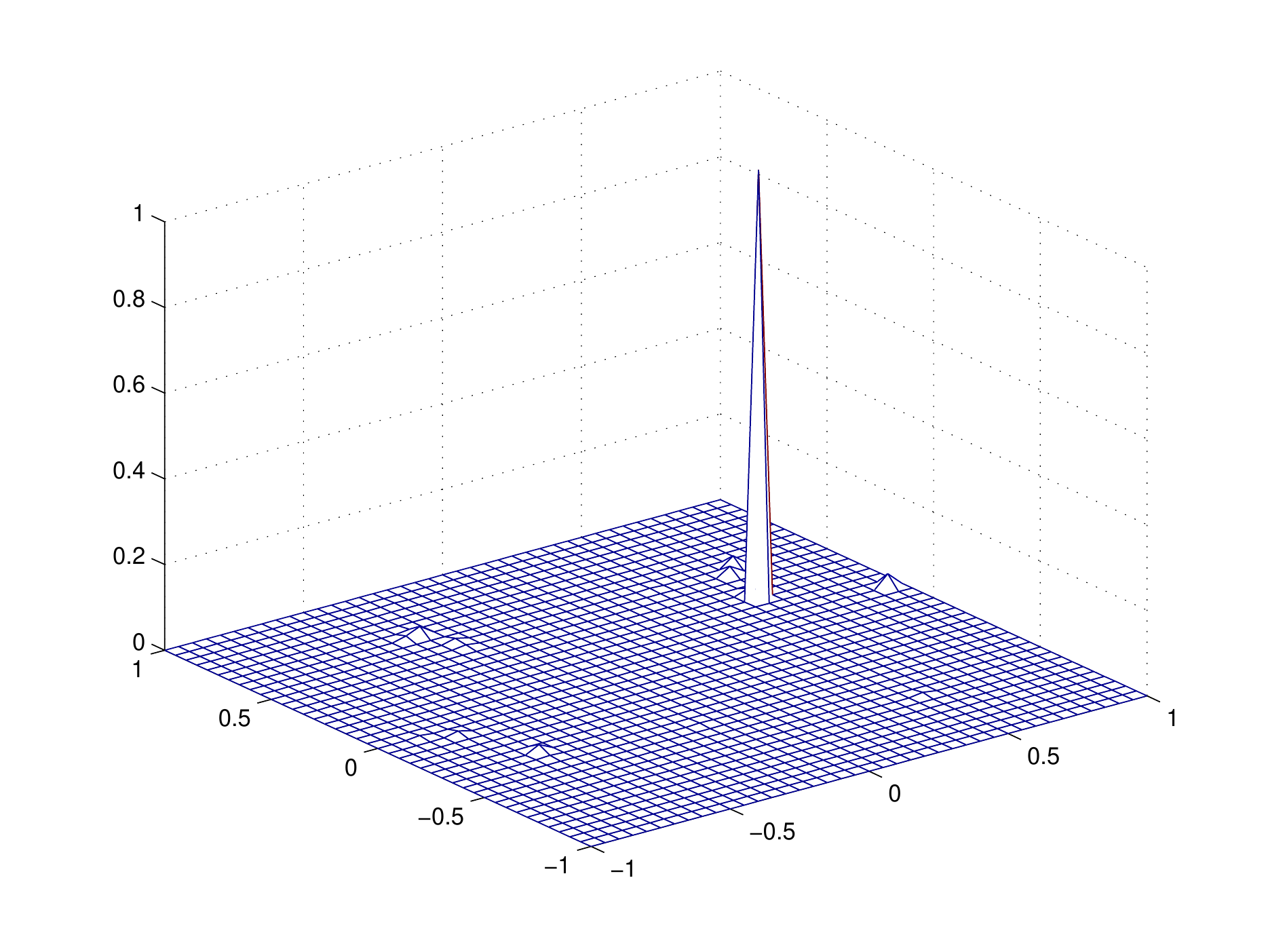}}
\end{indented}
\end{figure}

Standard implementation of the Gibbs sampler in this problem will not work,
since we discover, as Figure. \ref{fig:multimode} shows, that the posterior distributions are extremely multimodal.
The reason for this multimodality is clear. Indeed, what the algorithm essentially does is to look at concentrations of trajectories at different locations. If the threshold is set too low, as in the left part of Figure. \ref{fig:multimode}, one expects to find (and indeed finds) such concentrations in quite a few places.

To overcome this problem, after some experimentation we adopted a parallel tempering method
in order to improve mixing of simulations from this multimodal distribution.
%

We describe now the algorithm to sample from $f(\psi|\widetilde{Y}, M_2)$ in more details.
The reader interested only in the results of the implementation, can skip the following sub-sections and move directly to Section \ref{sec:simu}.

\subsubsection{Implementing the Parallel Tempering Algorithm} \label{subsec:pt}

One can find discussions of tempering algorithms
in \cite{GoLi} and \cite{Lia}.

We run $N$ parallel chains, each with equilibrium $f_i(x) \propto f(x)^{1/T_i}$, where $f(x)$ is the target posterior distribution $f(p, r, u | \widetilde{Y}, M_2)$ and $T_i$ is a given temperature level.
The temperature ladder $T_1 > \cdots > T_N=1$ plays the most important role in the algorithm, and is constructed in the following (trial-and-error) manner. We decide first the highest temperature such that a single MCMC run (e.g. using Metropolis-Hastings algorithm) at that temperature can explore the whole sample space easily (e.g. the acceptance rate of MH is about 90\%). Then the next temperature level is chosen such that the rate of exchanging samples with the chain at previous temperature is moderate (e.g. 20\%).
We have found in numerical experiments that $N=6$ works well, with highest temperature being 5 and exponentially decreasing to 1.
In our MCMC simulation, the Gibbs algorithm is implemented at each chain, and all chains start with random values. Let $\widetilde{x}^{(t)}=(x_1^{(t)}, \cdots, x_N^{(t)})$ denote the current population of samples from $N$ chains.

During each iteration of MCMC, the following steps of mutation and exchange are implemented.\\
\noindent\underline{Mutation}: Update $x_i^{(t)}$ to $x_i^{(t+1)}$ by the Gibbs sampler, for $i=1, \cdots, N$.
				Details are shown in the subsection \ref{subsec:gibbs}.\\
\noindent\underline{Exchange}: Starting from the first chain, try to swap with neighbors as follows,
        \begin{itemize}
            \item For $i=1, N$, exchange with the only neighbor; for $i=2, \cdots, N-1$, exchange with the two neighbors with equal probability.
            \item Then accept the exchange of states $i$ and $j$ with probability
            \begin{equation*}
                 \min \left\{  1, \exp\left(  \left[   \log f(x_j^{(t+1)}) - \log f(x_i^{(t+1)}) \right].
                \left[ T_i^{-1} - T_j^{-1}  \right] \right)  \right\}.
            \end{equation*}
        \end{itemize}
The chain with $T_N=1$, which has the target posterior distribution as equilibrium, is used in the harmonic mean estimate of Bayes factors.

\subsubsection{Implementing the Gibbs Sampler} \label{subsec:gibbs}

In each Gibbs update, the target distribution is one of the $f_i$'s. In the following context, a function $h(x|\cdot)$ refers to the full conditional distribution of $X$, given all the other unknown variables. Notice that the unknown additive constants in the logarithm of a distribution do not affect Gibbs sampler.


\noindent\underline{Joint Distribution in Each Iteration}:
It is easily seen that the joint posterior distribution of $(p,r,u)$ is given as
\bse
    \log \{f_i(p, r, u | \widetilde{Y} ) \}
         & = &  T_i^{-1}              \log \{ f( p, r , u | \widetilde{Y} ) \}   + C \\
         & = &  T_i^{-1} [ \log(r) + \log \{ f( \widetilde{Y} | p, r , u ) \} ]+ C.
\ese

\noindent\underline{Updating $p$}:
It is easily seen that
\bse
    \log \{ f_i(p| \cdot ) \} & = & T_i^{-1} \{ (n-J)\log (1-p) + J \log( pd^{-1}+1-p )    \}   + C_p.
\ese
To update $p$ for each $f_i(p|\cdot)$, the following steps are taken.
\begin{itemize}
    \item Compute the full conditional distribution $f_i(p_j|\cdot)=\pi_j$, for $j=1,\cdots,h $.
    \item Form $\omega_j = {\pi_j} / {\hbox{$\sum_{k=1}^h$} \pi_k}$, for $j=1,\cdots,h $.
    \item Draw from the vector $(p_1, \cdots, p_h)$ with probabilities $(\omega_1, \cdots, \omega_h)$.
\end{itemize}

\noindent\underline{Updating the Radius Component $r$ in the Polar Coordinates}:
It is easily seen that
\bse
    \log \{ f_i(r| \cdot ) \} & = & T_i^{-1} \{ J \log( pd^{-1}+1-p) - J\log (1-p) + \log(r) \} + C_r.
\ese
To update $r$, the Metropolis-Hastings algorithm is implemented. The proposal is a normal distribution $\mbox{N}(r_{\rm curr}$, $\sigma_{\rm prop, r}$) truncated on the interval $[0, 1]$, where the mean $r_{\rm curr}$ is the current value and the standard deviation $\sigma_{\rm prop, r}$ is a constant.

\noindent\underline{Updating the Angle Component $u$ in the Polar Coordinates}:
It is easily seen that
\bse
    \log \{ f_i(u| \cdot ) \} & = & T_i^{-1} J \{ \log( pd^{-1}+1-p) - \log (1-p) \} + C_u.
\ese
To update $u$, again the Metropolis-Hastings algorithm is implemented. Notice that the target function is periodic in $u$, and $u$ is restricted to $[0, 2\pi]$. To ease movement across boundary, we propose a value $u_{\rm prop}$ from the normal distribution $\mbox{N}(u_{\rm curr}, \sigma_{\rm prop, u}^2)$, where the mean $u_{\rm curr}$ is the current value and standard deviation $\sigma_{\rm prop, u}$ is a constant.
If $u_{\rm prop} > 2\pi$ or $u_{\rm prop} <0$, then the candidate is reset to be $u_{\rm prop, actual}=u_{\rm prop}\mbox{mod}(2\pi)$.

\subsection{Algorithm Summary}
Our Bayesian approach could be summarized as follows. Given a particular dataset $\widetilde{Y}$,  denote the posterior by $ f(\psi | \widetilde{Y}, M_2) $ as in equation (\ref{eq:post}),
\begin{enumerate}
 \item decide the hyperparameters $a_p$ and $b_p$ ($b_p>a_p>0$) based on prior knowledge of $p$, and the grid $a_p = p_1<\cdots<p_h = b_p$ according to desired precision in $p$,
 \item decide the temperature ladder $T_1 > \cdots > T_N=1$ and define
 $ f_i(\psi)=\{ f(\psi | \widetilde{Y}, M_2) \}^{1/T_i}$,
for $i=1, \cdots, N$,
 \item assign random initial values $\psi_{i,0} = (p_{i0}, r_{i0}, u_{i0})$ to each MCMC chain; set $t=0$,
 \item update $\psi_{i,t}$ to $\psi_{i,t+1}$ by Gibbs sampler as explained in subsection \ref{subsec:gibbs}, and exchange $\psi_{i,t+1}$ with its neighbor(s) as explained in subsection \ref{subsec:pt}, for $i=1, \cdots, N$; set $t=t+1$,
 \item repeat the last step until $t=K$; check the convergence and mixing of $\{\psi_{N, t}\}_{t=0}^{K}$; adjust the temperature ladder and repeat the above steps until the MCMC chain converges and mixes well,
 \item discard the first 20\% of the sequence $\{ \psi_{N, t} \}_{t=0}^{K}$ (called burn-in), take every 10 samples from the rest of the chain (called thinning), denote the new sequence by $\{ \psi_{j} \}$, which are samples from the posterior distribution,
 \item calculate the Bayes factor estimator $\widehat{\mbox{BF}}$ as in equation (\ref{eq:BFcalc}), and the posterior sample mean $\widehat{\psi} = (\widehat{p}, \widehat{r}, \widehat{u})$,
 \item conclude the presence of the source if $\widehat{\mbox{BF}}>3$; the source strength (i.e. the emission rate) is estimated by $\widehat{p}$; if $\widehat{\mbox{BF}}>3$, the location of the detected source is estimated by $\{ \widehat{r}\cos (\widehat{u}), \widehat{r}\sin(\widehat{u}) \}$.
\end{enumerate}
Furthermore, the uncertainty in estimation could be summarized by other statistics such as sample standard deviation.

\section{Simulation Study} \label{sec:simu}

We considered the situation where the size of the possible source is approximately known and is small compared to the size of the whole object. After choosing appropriate units, we assume that the object is the unit disk. The practically reasonable assumption is that the source radius is around $1\%$ of the object radius, i.e. $d=0.01$ (e.g., the object has dimension of several meters, while the source is of diameter of a few centimeters). The simulation is designed to examine the performance of the method at various emission rate levels, which are chosen as $p= 0.01$, $p = 0.005$, $p = 0.001$ and the case that no source exists, $p=0.00$. We experimented with the number of detected particles being $n=2\times10^5$ and $n = 5\times10^5$. With a fixed $d$, as the true emission rate $p$ (and thus signal-to-noise ratio) decreases, a larger total number of all detected particles is required for detection of a source. This can be explained by a simple application of the Central Limit Theorem CLT (see, e.g. \cite{ADHK}). The prior values $[0<p_1, \dots, p_h]$ are assumed to be located near the true value of $p$, which in many applications is known with some uncertainty. At each level $p$, $10$ simulated data sets where generated and analyzed, including also the case $p=0$. The results with two different sample sizes are summarized in Table \ref{table:levels} and Table \ref{table:largen}, respectively. Along with the Bayes factors, we also report the posterior probability that there is no sources, namely
\bse
\pr (p=0|\widetilde{Y}) &=& \pr (M_1|\widetilde{Y}) \\
                        &=& \frac{\pr(M_1)\pr(\widetilde{Y}|M_1)}
                            {\pr(M_1)\pr(\widetilde{Y}|M_1)+\pr(M_2)\pr(\widetilde{Y}|M_2)}
                         = (1+{\rm BF})^{-1},
\ese
where BF refers to the Bayes factor.

\begin{table}[t]
\caption{\label{table:levels}  
Summary of Bayes factors for simulation in Section \ref{sec:simu}. Sample size $n=2\times10^5$. There are 10 simulations performed at each combination of level $p$ and location, and 20 simulations at $p=0$. The values reported in the table are minimum, median, maximum of 10 Bayes factors, and the proportion of Bayes factors being greater than 3. In the last column is the median of $\pr(p=0|\widetilde{Y})$ calculated from the 10 data sets.}
\begin{indented}
\lineup
\item[] \begin{tabular}{@{}lllllll}
\br
$p$       &  Location     &  Min                    &  Med                   &     Max                 & Prop$>$3  &  $\pr(p=0|\widetilde{Y})$
\\
\mr
0.01      &  (0.6, 0.3)   &  3.7$\times10^{282}$ &  Inf                   &     Inf                 &    1      & 0 \\
0.01      &  (0.96, -0.1) &  6.7$\times10^{297}$ &  Inf                   &     Inf                 &    1      & 0 \\
0.005     &  (0.6, 0.3)   &  2.8$\times10^{75}$  &  3.1$\times10^{89}$ &  2.4$\times10^{99}$  &    1      & 3.1$\times10^{-90}$ \\
0.005     &  (0.96, -0.1) &  3.5$\times10^{76}$  &  1.2$\times10^{86}$ &  4.6$\times10^{106}$ &    1      & 8.1$\times10^{-87}$ \\
0.001     &  (0.6, 0.3)   &  0.3299                 &  11.2              &  7.7$\times10^{5}$   &   0.6     & 8.1$\times10^{-2}$ \\
0.001     &  (0.96, -0.1) &  0.1032                 &  1.4                &  2.1$\times10^{5}$   &   0.4     & 0.4051 \\
0         &   n/a         &  0.36                 &  1.06               &  24.73                  &   0.1     & 0.48 \\
\br
\end{tabular}
\end{indented}
\end{table}

One can note from the results that if $p=0.005$ or $p=0.01$, much smaller sample sizes are sufficient to detect the source. As the level decreases, say for $p=0.0005$, much larger sample sizes are required. In particular,  the rows with $p=0.001$ of the Table \ref{table:levels} show that sensitivity is not too high. The reason is that the number of detected particles, $n=2\times10^5$, is not high enough. The next table shows significant improvements when the number of particles is increased.

\begin{figure}[t]
\caption{ \label{fig:estloc} 
The estimated location of the source with $n=5\times10^{5}$ sample counts is shown for various emitting levels, with 95\% highest posterior density region. 
The top left plot is with $p=0.01$; the top right plot is with $p=0.005$;
the bottom left plot is with $p=0.001$; the bottom right plot is with $p=0$.
In the plots where a source exists, the location of its center is indicated by the intersection of gray dashed lines.
Each of the above figures is plotted with the posterior sample having median Bayes factors among the 10 simulated cases.
}
\begin{indented}
\item[]
  \begin{minipage}[t]{0.2\textwidth}
    \includegraphics[width=2in]{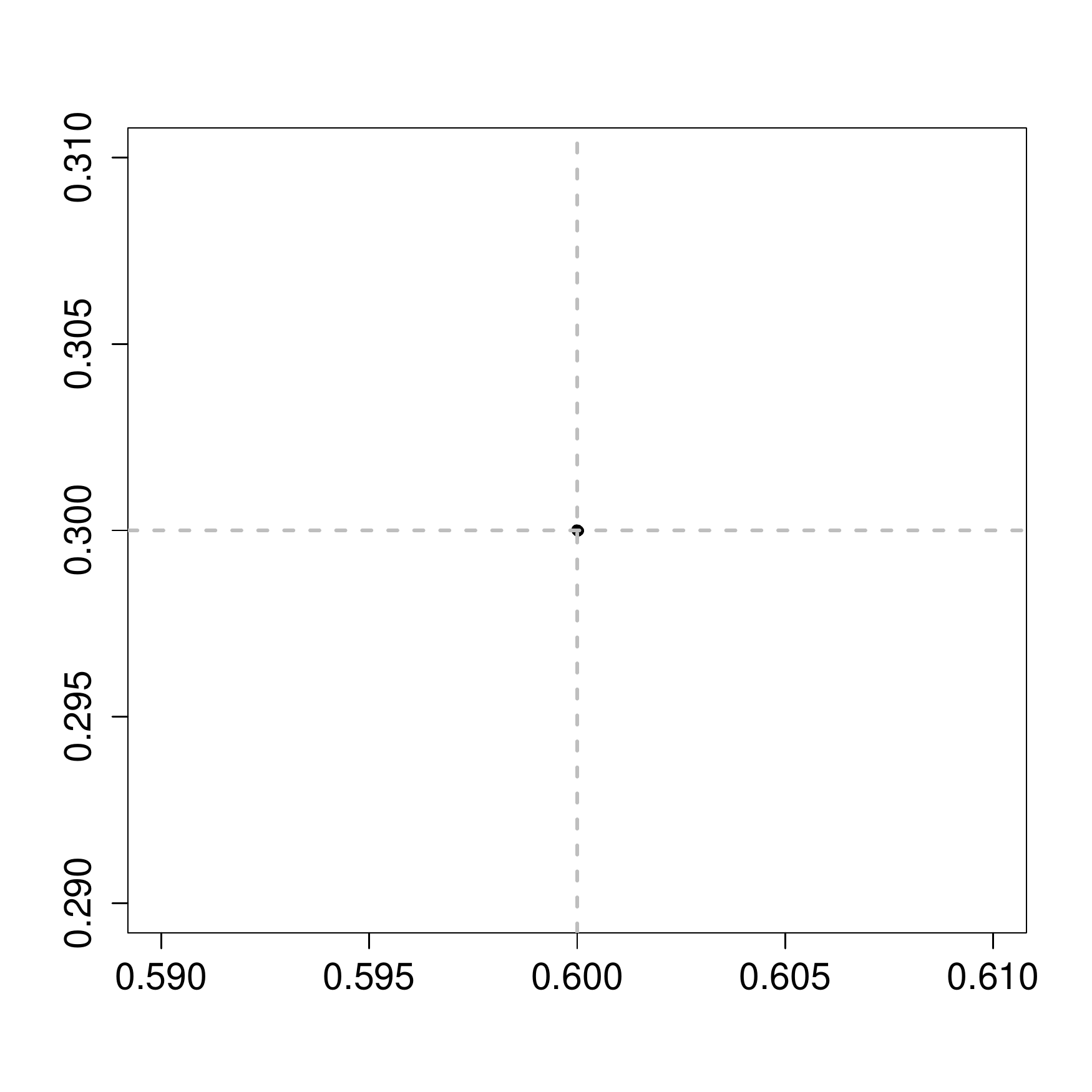}
  \end{minipage}
  \hspace{0.2\textwidth}
  \begin{minipage}[t]{0.2\textwidth}
    \includegraphics[width=2in]{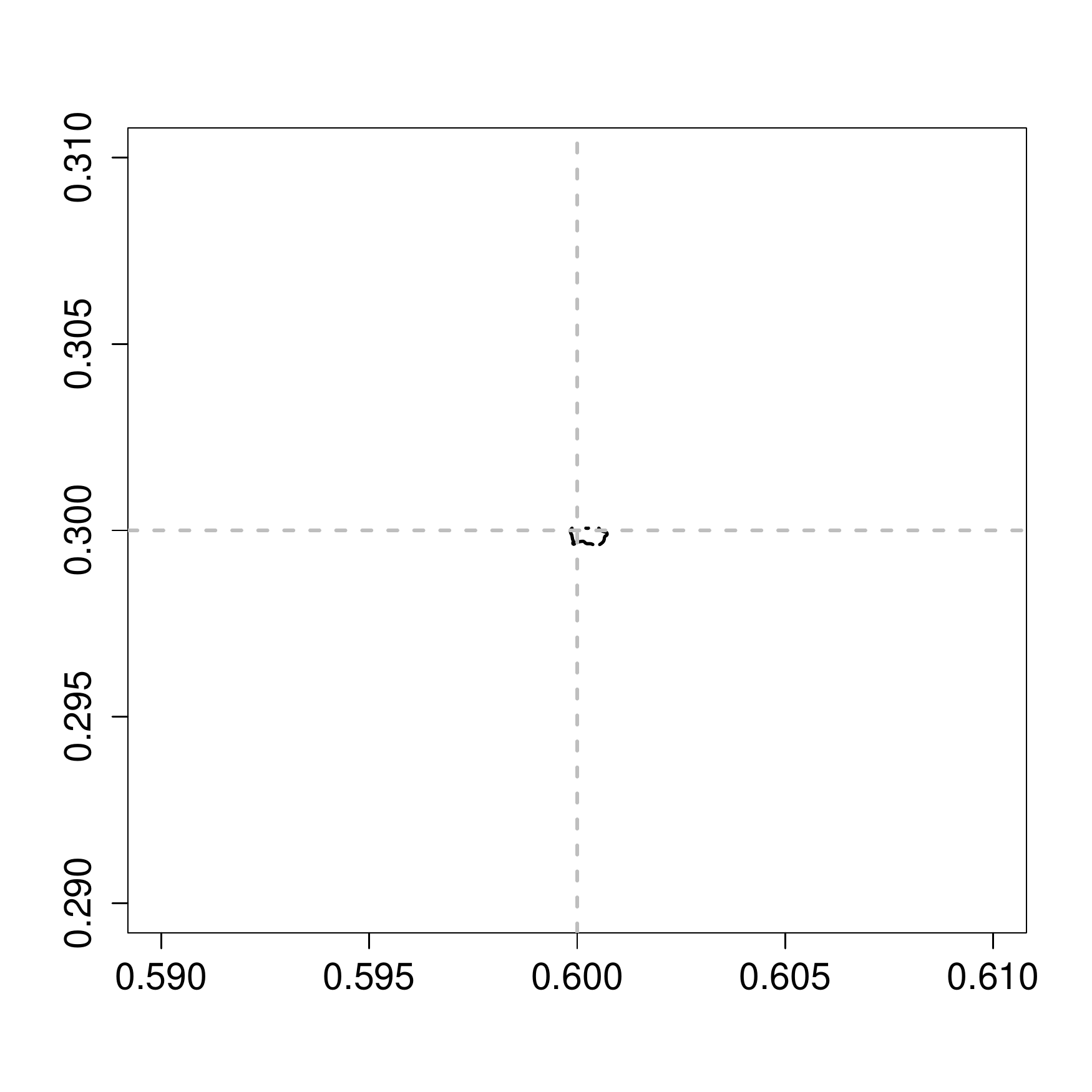}
  \end{minipage}                               \\[2 pt]
  \begin{minipage}[b]{0.2\textwidth}
    \centering
    \includegraphics[width=2in]{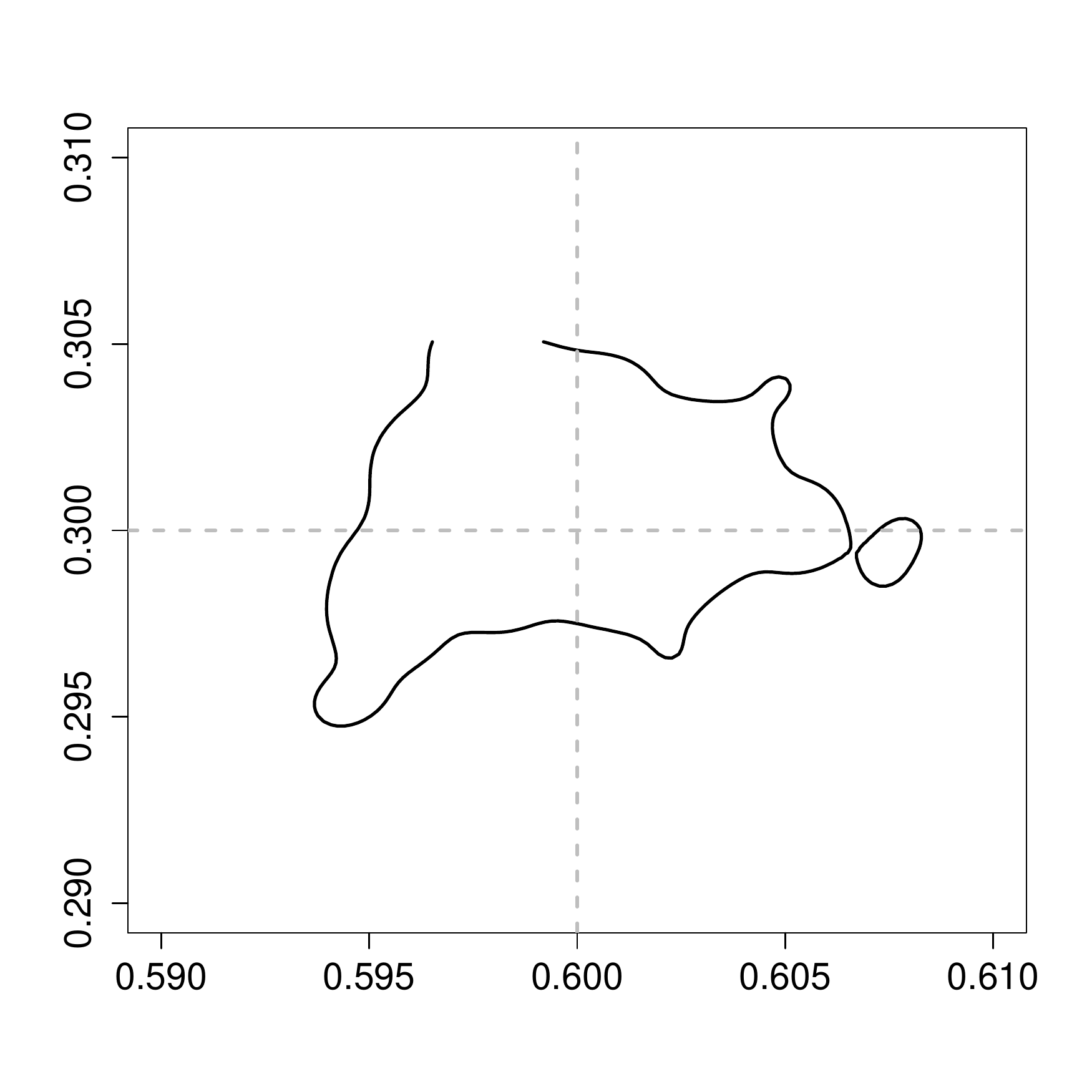}
  \end{minipage}
  \hspace{0.2\textwidth}
  \begin{minipage}[b]{0.2\textwidth}
    \centering
    \includegraphics[width=2in]{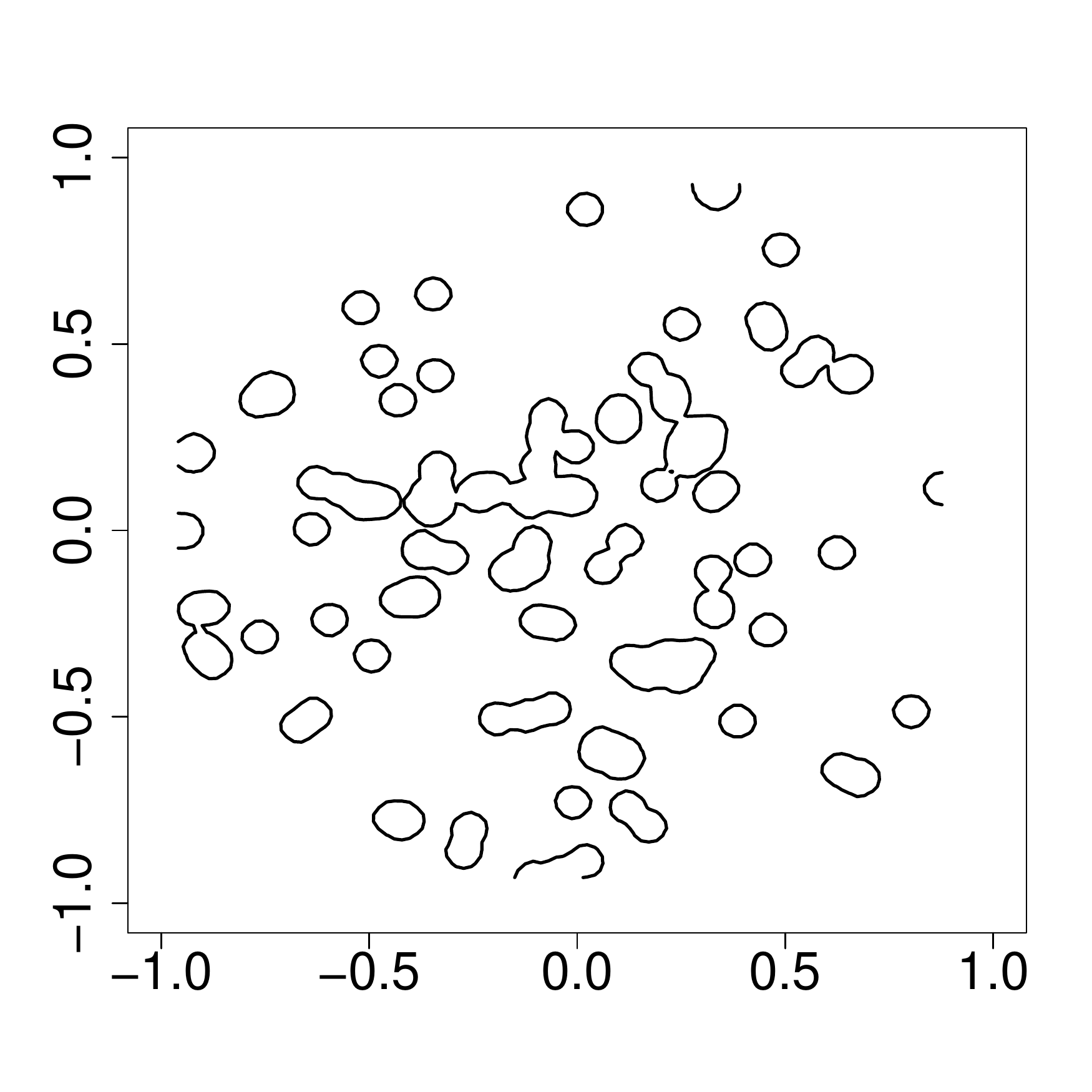}	
  \end{minipage}%
\end{indented}
\end{figure}

The results in Table \ref{table:largen} clearly show very high sensitivity to the presence of a source, which is indicated by the overall large values of BF. Furthermore, the location, if the source is present, can be also found with high accuracy, as shown in Figure \ref{fig:estloc}.
\begin{table}[t]
\caption{\label{table:largen} Repeat of the Table \ref{table:levels} with $5\times 10^{5}$ samples.}
\begin{indented}
\lineup
\item[] \begin{tabular}{@{}lllllll}
\br
$p$     &  Location     &  Min     		&  Med                  &   Max               	& Prop$>$3  &  $\pr(p=0|\widetilde{Y})$  \\
\mr
0.01	& (0.6,0.3)	& Inf			& Inf			 & Inf			 & 1	&	0	\\
0.01	& (0.96, -0.1)	& Inf			& Inf			 & Inf			& 1	&	0	\\
0.005	& (0.6,0.3)	& 1.02$\times 10^{219}$	& 9.81$\times 10^{235}$	& 7.77$\times 10^{275}$	& 1	&	6.48$\times 10^{-236}$\\
0.005	& (0.96, -0.1)	& 2.00$\times 10^{214}$	& 4.13$\times 10^{223}$	& 1.91$\times 10^{255}$	& 1	&	 3.35$\times 10^{-223}$\\
0.001	& (0.6,0.3)	& 523.76		& 1.41$\times 10^{9}$	 & 5.44$\times 10^{16}$	& 1	&	 7.32$\times 10^{-10}$\\
0.001	& (0.96, -0.1)	& 41.68 		& 2.35$\times 10^{9}$	& 1.60$\times 10^{12}$	& 1	&	 4.96$\times 10^{-10}$\\
0	& n/a	        & 0.65			& 1.16			 & 2.23			 & 0	& 0.46	\\
\br
\end{tabular}
\end{indented}
\end{table}

\section{Concluding Remarks}
\begin{itemize}
 \item The results of this study show that Bayesian methods can be successfully used for detection of low emission small sources in the cases of realistic parameters. High sensitivity can be achieved if the observation time (and thus total count of particles detected) is sufficiently large.

\item The assumption that the detector can determine the directional information is crucial. We believe that otherwise detection with such low values of SNR (signal-to-noise ratio) would be impossible. The determination of the incoming direction is usually achieved by detector collimation. This is not an option with the extremely low SNR levels, since it would most probably eliminate completely the useful signal. However, there exist the so called Compton type cameras for detecting $\gamma$-photons (e.g., \cite{ADHK,Comp} and references therein), as well as their analogs (although based upon a somewhat different physics) for neutron detection \cite{shield,SpCharlton}. These cameras do not use collimation, but can determine some less precise directional information. Namely, the camera is able to provide a (hollow) cone of possible directions of the incoming particle. Although this is a less precise (and highly over-determined) information, it is known (e.g., see \cite{ADHK} and references therein) how to convert the Compton type data into the precise directional information. Thus, the algorithms described can be used in conjunction with Compton type detectors.

\item It is planned to address in the future study the effects of scattering, inhomogeneous random background noise, Compton type cameras, as well as the $3D$ situation. It is also planned to compare the results and the computational cost of the Bayesian approach with the more analytic techniques of \cite{ADHK}.

\end{itemize}
\section*{Acknowledgments}

Carroll's research was supported by a grant from the National Cancer Institute (R37-CA057030).
Mallick's research was supported by a grant from the NSF DMS grants 0914951.
Work of Kuchment was partly supported by the DHS Grant 2008-DN-077-ARI018-04 and by the NSF DMS grants 0604778 and 0908208.
Work of all authors was supported in part by the IAMCS.
Thanks also go to W.~Charlton, Y.~Hristova, G.~Spence, and the reviewers for useful discussions and comments.

\end{document}